\documentclass[12pt,reqno]{amsart}
\usepackage{fullpage}
\usepackage{times}

\newtheorem{theorem}{Theorem} %[section]

\renewcommand{\mod}[1]{{\ifmmode\text{\rm\ (mod~$#1$)}\else\discretionary{}{}{\hbox{ }}\rm(mod~$#1$)\fi}}

\newsymbol\nmid 232D

\newcommand{\lcm}{\mathop{\rm lcm}}

\begin{document}

\title{A product of Gamma function values at fractions with the same denominator}
%\author[Greg Martin]{Greg Martin \\ Draft: \today}
\author{Greg Martin}
\address{Department of Mathematics \\ University of British Columbia \\ Room 121, 1984 Mathematics Road \\ Canada V6T 1Z2}
\email{gerg@math.ubc.ca}
\maketitle

\begin{abstract}
We give an exact formula for the product of the values of Euler's Gamma function evaluated at all rational numbers between 0 and 1 with the same denominator in lowest terms; the answer depends on whether or not that denominator is a prime power. A consequence is a surprisingly nice formula for the product of value of the Gamma function evaluated at the points of a Farey sequence.
\end{abstract}

{\em {\bf Note:} Since writing this note, I have been informed that Theorem~\ref{main theorem} was already proved by S\'andor and T\'oth~\cite{scooped}.}

The purpose of this note is to establish the following classical-seeming theorem concerning Euler's $\Gamma$-function evaluated at fractions that have the same denominator in lowest terms. The statement of the theorem uses (coincidentally) Euler's function $\phi(n)$, the number of integers between 1 and $n$ that are relatively prime to $n$, as well as von Mangoldt's function $\Lambda(n)$, defined to be $\log p$ if $n=p^r$ is a prime or a power of a prime and 0 otherwise.

\begin{theorem}
\label{main theorem}
For $n\ge2$,
\[
\prod_{\substack{k=1 \\ (k,n)=1}}^n \Gamma\big( \tfrac kn \big) = \frac{(2\pi)^{\phi(n)/2}}{\exp(\Lambda(n)/2)} = \begin{cases}
(2\pi)^{\phi(n)/2}/\sqrt p, &\text{if $n=p^r$ is a prime power}, \\
(2\pi)^{\phi(n)/2}, &\text{otherwise.}
\end{cases}
\]
\end{theorem}

A few special cases of this theorem have been noted before ($n=2,3,4,6$ for example), and it follows for prime $n$ from equation~\eqref {mult formula special case} below. Nijenhuis~\cite[page 4]{nijenhuis} established, by a more indirect method, the special case of the theorem where $n\equiv2\mod4$.

\begin{proof}
Gauss's multiplication formula~\cite[equation~(3.10)]{artin} says that
\[
\prod_{k=0}^{n-1} \Gamma\big( \tfrac{z + k}n \big) = (2\pi)^{(n-1)/2} n^{1/2-z} \Gamma(z)
\]
for any complex number $z$ for which both sides are defined; taking $z=1$ yields
\begin{equation}
\prod_{k=1}^n \Gamma\big( \tfrac kn \big) = (2\pi)^{(n-1)/2} n^{-1/2}.
\label{mult formula special case}
\end{equation}
Define the two functions
\[
F(n) = \sum_{k=1}^n \log\Gamma\big( \tfrac kn \big) \quad\text{and}\quad R(n) = \sum_{\substack{k=1 \\ (k,n)=1}}^n \log\Gamma\big( \tfrac kn \big).
\]
It is immediate from these definitions that $F(n) = \sum_{d\mid n} R(d)$; hence M\"obius inversion~\cite[second displayed equation after equation~(2.10)]{MV} yields
\[
R(n) = \sum_{d\mid n} \mu(d) F\big( \tfrac nd \big).
\]
From equation~\eqref{mult formula special case} we see that $F(n) = \log \big( (2\pi)^{(n-1)/2} n^{-1/2} \big)$, and so
\begin{align*}
R(n) &= \sum_{d\mid n} \mu(d) \bigg( \frac{n/d-1}2 \log 2\pi - \tfrac12\log \tfrac nd \bigg) \\
&= \tfrac12\log2\pi \sum_{d\mid n} \mu(d) \tfrac nd - \tfrac12\log2\pi \sum_{d\mid n} \mu(d) - \tfrac12 \sum_{d\mid n} \mu(d) \log \tfrac nd.
\end{align*}
Each of these three divisor sums is standard in number theory (see~\cite{MV}, where they appear as the first displayed equation in the proof of Theorem 2.1, equation~(1.20), and the displayed equation before equation~(2.10), respectively): as long as $n\ge2$, we have
\[
R(n) = \big( \tfrac12\log2\pi \big) \phi(n) - 0 - \tfrac12 \Lambda(n).
\]
Taking exponentials of both sides establishes the theorem.
\end{proof}

It was known in the nineteenth century that the geometric mean of the $\Gamma$ function on the interval $(0,1]$ is $\sqrt{2\pi}$, in the sense that
\[
\int_0^1 \log\Gamma(x)\,dx = \tfrac12\log2\pi.
\]
(One can deduce this, for example, by integrating the Weierstrass formula~\cite[equation~(2.9)]{artin}
\[
\log\Gamma(z) = -\gamma z-\log z + \sum_{j=1}^\infty \big( \tfrac zj - \log\big( 1+\tfrac zj \big) \big)
\]
term by term; another proof uses the reflection formula $\Gamma(z)\Gamma(1-z) = \pi\csc\pi z$ together with a known evaluation of the integral $\int_0^{1/2} \log(\sin\pi x)\,dx$.) Therefore if we multiply together $n$ values of the $\Gamma$ function on points in this interval, we would expect the product to be comparable to $(2\pi)^{n/2}$. We can deduce from first principles that the product will be less than $(2\pi)^{n/2}$ if we sample the $\Gamma$ function at $\frac1n, \frac2n, \dots, \frac nn$, since $\Gamma$ is decreasing on that interval; in fact, equation~\eqref {mult formula special case} tells us that the product will be less by a factor of precisely $1/\sqrt{2\pi n}$. Applying equation~\eqref {mult formula special case} twice, at $2n$ and $n$, and dividing shows that we do better to sample at the midpoints, rather than the right-hand endpoints, of $n$ intervals of equal length:
\begin{equation}
\prod_{k=1}^n \Gamma\big( \tfrac{2k-1}{2n} \big) = \frac{(2\pi)^{n/2}}{\sqrt2}.
\end{equation}
Theorem~\ref{main theorem} tells us that sampling at the $\phi(n)$ points $\{ \frac kn\colon 1\le k\le n,\, (k,n)=1\}$ curiously gives us exactly the default expectation $(2\pi)^{\phi(n)/2}$, unless $n$ is a prime power.

Finally, we comment that the $\Lambda$-function satisfies the identity~\cite[Section 2.2.1, exercise 1(a)]{MV}
\begin{equation*}
\sum_{n=1}^N \Lambda(n) = \log \big( \lcm[1, 2, \dots, N] \big).
\end{equation*}
This allows us to compute the product of the $\Gamma$-function sampled over points in a Farey sequence. Let $F_N$ denote the set of all rational numbers in the open interval $(0,1)$ whose denominator in lowest terms is at most $N$ (note that usually one includes the fractions $\frac01$ and $\frac11$ in this Farey sequence, but here we do not). Applying Theorem~\ref{main theorem} to $n=2,3,\dots,N$ and multiplying the identities together yields the formula
\begin{equation}
\prod_{r\in F_N} \frac{\Gamma(r)}{\sqrt{2\pi}} = \big( \lcm[1, 2, \dots, N] \big)^{-1/2}.
\end{equation}
It was noted by Luschny and Wehmeier~\cite{LW} that this last equation is equivalent, via the reflection formula $\Gamma(z)\Gamma(1-z) = \pi\csc\pi z$, to the identity
\[
\lcm[1, 2, \dots, N] = \tfrac12 \bigg( \prod_{\substack{r\in F_N \\ r \le 1/2}} 2\sin\pi r \bigg)^2;
\]
in fact they found an alternate proof using cyclotomic polynomials.

\end{document}